\documentclass[review]{elsarticle}
\usepackage{graphicx}
\usepackage{lineno,hyperref}
\usepackage{subfigure}
\usepackage{amsmath, amssymb}
\usepackage{amsbsy}
\usepackage{bm}
\usepackage{makecell}
\usepackage{color}
\usepackage{float}
\usepackage{multirow}
\usepackage{mathtools}
\usepackage{amsmath,amssymb,amsfonts,epsfig,subfigure}
\usepackage{xcolor}
\graphicspath{{fig/new/}}
\usepackage[demo,abs]{overpic}
\definecolor{newcolor}{rgb}{.8,.349,.1}
\usepackage[linesnumbered,lined,boxed,commentsnumbered]{algorithm2e}

\modulolinenumbers[5]
\usepackage[margin=1. in]{geometry}
\journal{Elsevier}

\begin{document}

\begin{frontmatter}

\title{Development of A Hermite Weighted Compact Nonlinear Scheme based on the Two-Stage Fourth-Order Temporal Accurate Framework\\
{Report \#1} }


\date{\today}
\author[mymainaddress]{Huaibao Zhang}
\ead{zhanghb28@sysu.edu.cn}

\author[mysecondaryaddress]{GX W\corref{mycorrespondingauthor}}
\ead{xuchg5@sysu.edu.cn}
%
%
%
\address[mymainaddress]{School of Aeronautics and Astronautics, Sun Yat-sen University, Guangzhou 510006, China }
%
%
\begin{abstract}
Improved five-point low dissipation nonlinear schemes are proposed in this paper within the framework of weighted compact nonlinear schemes (WCNSs)~\cite{Deng2000}. Particularly we follow the work of Li and Du~\cite{Li2016} on the two-stage fourth-order temporal accurate discretization scheme, which is developed based on the Lax-Wendroff method.

\end{abstract}


\end{frontmatter}

\linenumbers

\section{Introduction} \label{sec:Intro}

\section{Fundamentals of the numerical methods} \label{sec:Fund}

Consider the time-dependent hyperbolic conservation law in one dimension together with its initial condition, given by

\begin{equation}\label{eqn:gov1}
	\begin{aligned}
		&\frac{\partial \mathbf{u}}{\partial t}+\frac{\partial \mathbf{f}(\mathbf{u})}{\partial x} +\frac{\partial \mathbf{g}(\mathbf{u})}{\partial y}=0, \quad (x, y)  \in \mathbb{R} \times \mathbb{R}, \quad t>0 \\
		&\mathbf{u}(x, 0)=\mathbf{u}_{0}(x), \quad x \in \mathbb{R},
	\end{aligned}
\end{equation}
where $\mathbf{u}$ denotes the vector of conservative variables, and $\mathbf{f}(\mathbf{u})$ the vector of flux terms.

\begin{equation}\label{eqn:govxx}
\begin{aligned}
	\left(\frac{\partial \mathbf{u}}{\partial t}\right)_{i+1 / 2, j}=&-\mathbf{A}^{+}\left(\mathbf{u}_{i+1 / 2, j}\right)\left(\frac{\partial \mathbf{u}}{\partial x}\right)_{i+1 / 2, j}^{L}-R I^{+} L\left(\frac{\partial \mathbf{g}\left(\mathbf{u}_{i+1 / 2, j}^{L}\right)}{\partial y}\right)_{i+1 / 2, j} \\
	&-\mathbf{A}^{-}\left(\mathbf{u}_{i+1 / 2, j}\right)\left(\frac{\partial \mathbf{u}}{\partial x}\right)_{i+1 / 2, j}^{R}-R I^{-} L\left(\frac{\partial \mathbf{g}\left(\mathbf{u}_{i+1 / 2, j}^{R}\right)}{\partial y}\right)_{i+1 / 2, j}
\end{aligned}
\end{equation}

\begin{equation}
	\Omega_{i j} \equiv\left(x_{i-1 / 2}, x_{i+1 / 2}\right) \times\left(y_{i-1 / 2}, y_{i+1 / 2}\right)
\end{equation}

\begin{equation}
	\mathbf{f}(\mathbf{u})=\frac{1}{\Delta x} \int_{x-\Delta x / 2}^{x+\Delta x / 2} \mathbf{h}(\xi,y, t) d \xi
\end{equation}

\begin{equation}
	\mathbf{g}(\mathbf{u})=\frac{1}{\Delta y} \int_{y-\Delta y / 2}^{y+\Delta y / 2} \mathbf{h}(x,\xi, t) d \xi
\end{equation}

\begin{equation}
	\frac{d}{d t} \mathbf{u}_{i}(t)=\Re\left(\mathbf{u}_{i}\right)=-\frac{\mathbf{h}_{ i+1 / 2,j}-\mathbf{h}_{i-1 / 2,j}}{\Delta x}-\frac{\mathbf{h}_{i,j+1 / 2}-\mathbf{h}_{i,j-1 / 2}}{\Delta y}
\end{equation}

\begin{equation}
	\frac{d}{d t}\left(\frac{d}{d t} \mathbf{u}_{i}(t)\right)=\frac{\partial}{\partial t} \Re\left(\mathbf{u}_{i}\right)=-\frac{\mathbf{h}_{i+1 / 2,j}^{\prime}-\mathbf{h}_{i-1 / 2,j}^{\prime}}{\Delta x} -\frac{\mathbf{h}_{i,j+1 / 2}^{\prime}-\mathbf{h}_{i,j-1 / 2}^{\prime}}{\Delta y}
\end{equation}

\begin{equation}
	\begin{aligned}
		\mathbf{h}_{i+1 / 2, j}^{\prime} \approx & \frac{3}{640}\left(\frac{\partial \mathbf{f}}{\partial t}\right)_{i-3 / 2,j}-\frac{29}{480}\left(\frac{\partial \mathbf{f}}{\partial t}\right)_{i-1 / 2,j}+\frac{1067}{960}\left(\frac{\partial \mathbf{f}}{\partial t}\right)_{i+1 / 2,j} \\
		&-\frac{29}{480}\left(\frac{\partial \mathbf{f}}{\partial t}\right)_{i+3 / 2,j}+\frac{3}{640}\left(\frac{\partial \mathbf{f}}{\partial t}\right)_{i+5 / 2,j}
	\end{aligned}
\end{equation}

\begin{equation}
	\left(\frac{\partial \mathbf{f}}{\partial t}\right)_{i+1 / 2,j}=\left(\frac{\partial \mathbf{f}}{\partial \mathbf{u}}\left(\mathbf{u}_{i+1 / 2,j}\right)\right)\left(\frac{\partial \mathbf{u}}{\partial t}\right)_{i+1 / 2,j}
\end{equation}

Spatial discretization of Eq.~\eqref{eqn:gov1} is performed on an equally-spaced grid with distance between two adjacent grid nodes denoted by $h$.  At each node $i$, we define $x_i=ih, \; i = 1, \cdots, N$, and $\mathbf{u}_i=\mathbf{u}(x_i,t)$. The midpoints associated to the node $i$ are defined as $x_{i\pm\frac{1}{2}} = x_i \pm \frac{h}{2} $, which indicate the cell interfaces across which the fluxes are evaluated. 

Evaluation of Eq.~\eqref{eqn:gov1} at each node yields a system of ordinary differentiation equations, given in the form of 
\begin{equation}\label{eqn:semi-dis}
	\frac{d {\mathbf{u}}_{i}(t)}{d t}=\mathcal{R}_{i}(\mathbf{u}) = -\frac{\partial \mathbf{f}(\mathbf{u})}{\partial x} \vert_{x=x_i} \quad i = 1, \cdots, N.	
\end{equation}
where we define $\mathcal{R}_{i}(\mathbf{u})$ as the spacial operator.

A semi-discrete solution in the conservative
finite difference form can be further obtained as
\begin{equation}\label{eqn:semi-dis2}
	\frac{d {\mathbf{u}}_{i}(t)}{d t}=-\frac{\mathbf{h}_{i+1 / 2}-\mathbf{h}_{i-1 / 2}}{h},
\end{equation}
where the primitive function or the often so-called numerical flux function $\mathbf{h}(x,t)$ is implicitly defined by 
\begin{equation}\label{eqn:fluxfunch}
	\mathbf{f}(\mathbf{u})=\frac{1}{h } \int_{x-h / 2}^{x+h / 2} \mathbf{h}(\xi, t) d \xi . 
\end{equation}

In the framework of weighted compact nonlinear schemes (WCNSs) \cite{Deng2000}, the flux $\mathbf{h}_{i+1/2}$ is numerically approximated using the linear combination of mid-point flux terms. The approximated flux, denoted by $\widehat{\mathbf{h}}_{i+1/2}$, is evaluated by
\begin{equation}\label{eqn:linear}
	\widehat{\mathbf{h}}_{i+1/2} =\frac{3}{640} \widehat{\mathbf{f}}_{i-\frac{3}{2}}-\frac{29}{480} \widehat{\mathbf{f}}_{i-\frac{1}{2}}+\frac{1067}{960}\widehat{\mathbf{f}}_{i+\frac{1}{2}}-\frac{29}{480} \widehat{\mathbf{f}}_{i+\frac{3}{2}}+\frac{3}{640} \widehat{\mathbf{f}}_{i+\frac{5}{2}},  
\end{equation}
and Taylor series expansion of $\widehat{\mathbf{h}}_{i+1/2}$ indicates that 
\begin{equation}\label{eqn:taylor}
	\widehat{\mathbf{h}}_{i+1/2}=\mathbf{h}_{i+1/2}+O\left(h^{5}\right). 
\end{equation}

The unknown mid-point flux terms on the right hand side of Eq.~\eqref{eqn:linear} are computed using numerical upwind flux functions, which are given in a generic form 
\begin{equation} \label{eqn:flux}
	\widehat{\mathbf{f}}_{i+\frac{1}{2}}=\frac{1}{2}\left[\left({\mathbf{f}}(\mathbf{u}_{R,i+\frac{1}{2}})+{\mathbf{f}}(\mathbf{u}_{L,i+\frac{1}{2}})\right)-|\widehat{\mathcal{A}}_{i+\frac{1}{2}}|\left(\mathbf{u}_{R,i+\frac{1}{2}}-\mathbf{u}_{L,i+\frac{1}{2}}\right)\right]  , 
\end{equation}
where the high-order interpolated flow variables on the left- and right-hand sides of the mid-point $x_{i + \frac{1}{2}}$ are denoted by the subscripts, $L$ and $R$, respectively, and $\widehat{\mathcal{A}}_{i+\frac{1}{2}}$ denotes the approximate Jacobian matrix of the flux function with respect to the conservative variables, i.e., $\widehat{\mathcal{A}}_{i+\frac{1}{2}} \approx \frac{\partial \mathbf{f}}{\partial \mathbf{u}}\left(\mathbf{u}_{i+\frac{1}{2}}^{n}\right)$. 

Before proceeding to the section~\ref{sec:hermite} for the detailed discussion of nonlinear interpolation of aforementioned flow variables in Eq.~\eqref{eqn:flux}, we first review the time marching algorithm used to integrate Eq.~\eqref{eqn:semi-dis}, since it is the core building block of the numerical methods in the proposed work. Instead of using the traditional third-order strongly stable Runge--Kutta method~\cite{Gottlieb2001}, we follow the work of Li and Du~\cite{Li2016} on the two-stage fourth-order temporal accurate discretization scheme, which is developed based on the Lax-Wendroff method.

The basic idea in our reconstruction is to use only $r$ stencils to
reconstruct the point-wise values of solutions and spatial derivatives for the $2r-1$-order ADER scheme in one dimension, while in two dimensions, the dimension-by-dimension sub-cell reconstruction approach for spatial derivatives is employed

The hyperbolic conservation laws may develop discontinuities in its solution even if the initial
conditions are smooth. 

An important advantage or feature of the WCNS method is that the variable interpolation is performed on the primitive/conservative/characteristic variables from the solution points to the ux points.

\subsection{Two-stage fourth-order temporal accurate scheme}\label{sec:2s4f}

Given the flow state at time step $n$, the two-stage temporal scheme performed to integrate Eq.~\eqref{eqn:semi-dis} to the next time level $n+1$ is summarized as follows.

\begin{itemize}
	\item [{\bf Step 1.}] Define a vector of intermediate node values ${\mathbf{u}}_{i}^{n+\frac{1}{2}}$, which are calculated by
	\begin{equation}\label{eqn:inter1}
		\begin{aligned}
			&{\mathbf{u}}_{i}^{n+\frac{1}{2}}={\mathbf{u}}_{i}^{n}-\dfrac{k}{2 h}\left[\widehat{\mathbf{h}}_{i+\frac{1}{2}}^{*}-\widehat{\mathbf{h}}_{i-\frac{1}{2}}^{*}\right] \\
			&\widehat{\mathbf{h}}_{i+\frac{1}{2}}^{*}=\widehat{\mathbf{h}}_{i+\frac{1}{2}}^{n}+\dfrac{k}{4}  \left(\dfrac{\partial \widehat{\mathbf{h}}}{\partial t}\right)_{i+\frac{1}{2}}^{n} 
		\end{aligned}\quad .
	\end{equation}
	where $k$ is the time step size.
	
	\item [{\bf Step 2.}] Advance the solution to the time state $t^{n+1} = t^n +k$ by
	\begin{equation}\label{eqn:inter2}
		\begin{aligned}
			&{\mathbf{u}}_{i}^{n+1}={\mathbf{u}}_{i}^{n}-\dfrac{k}{ h}\left[\widehat{\mathbf{h}}_{i+\frac{1}{2}}^{4th}-\widehat{\mathbf{h}}_{i-\frac{1}{2}}^{4th}\right] \\
			&\widehat{\mathbf{h}}_{i+\frac{1}{2}}^{4th}=\widehat{\mathbf{h}}_{i+\frac{1}{2}}^{n}+\dfrac{k}{6}  \left(\dfrac{\partial \widehat{\mathbf{h}}}{\partial t}\right)_{i+\frac{1}{2}}^{n} + \dfrac{k}{3}\left(\dfrac{\partial \widehat{\mathbf{h}}}{\partial t}\right)_{i+\frac{1}{2}}^{n+\frac{1}{2}}
		\end{aligned}\quad .
	\end{equation}	
\end{itemize}

\vspace{2ex}
\noindent \textbf{Remark 1.} In {\bf Step 1}, the flux $\widehat{\mathbf{h}}_{i+\frac{1}{2}}^{n}$ can be readily determined from Eq.~\eqref{eqn:linear}, and its derivative with respect to time is also approximated using Eq.~\eqref{eqn:linear}, such that

\begin{equation}\label{eqn:dh2}
	\left(\dfrac{\partial \widehat{\mathbf{h}}}{\partial t}\right)_{i+\frac{1}{2}}^{n}  =\dfrac{3}{640} \left(\dfrac{\partial \widehat{\mathbf{f}}}{\partial t}\right)_{i-\frac{3}{2}}^{n} -\frac{29}{480} \left(\dfrac{\partial \widehat{\mathbf{f}}}{\partial t}\right)_{i-\frac{1}{2}}^{n} +\frac{1067}{960}\left(\dfrac{\partial \widehat{\mathbf{f}}}{\partial t}\right)_{i+\frac{1}{2}}^{n}-\frac{29}{480} \left(\dfrac{\partial \widehat{\mathbf{f}}}{\partial t}\right)_{i+\frac{3}{2}}^{n}+\frac{3}{640} \left(\dfrac{\partial \widehat{\mathbf{f}}}{\partial t}\right)_{i+\frac{5}{2}}^{n}. 
\end{equation}
The flux derivative $\left(\frac{\partial \widehat{\mathbf{h}}}{\partial t}\right)_{i+\frac{1}{2}}^{n+\frac{1}{2}}$ in {\bf Step 2} can be obtained in the same way but at the intermediate time level. We discuss the evaluation of the variable $\left(\frac{\partial \widehat{\mathbf{f}}}{\partial t}\right)_{i+\frac{1}{2}}^{n}$ in the following subsection.

\subsection{Generalized Riemann problem solver}\label{sec:GRP}
While the flux term $\widehat{\mathbf{f}}_{i+\frac{1}{2}}$ can be computed by a general Riemann solver, for instance, the Rusanov scheme~\cite{Rusanov1961}, or the low dissipation hybrid Rusanov-Roe scheme~\cite{Wang2016}, the evaluation of $\left(\frac{\partial \widehat{\mathbf{f}}}{\partial t}\right)_{i+\frac{1}{2}}$ plays a key role, which is conducted following the Cauchy-Kovalevskaya procedure, such that
\begin{equation}\label{eqn:ck}
\left(\frac{\partial \widehat{\mathbf{f}}}{\partial t}\right)_{i+\frac{1}{2}} = \left(\frac{\partial \widehat{\mathbf{f}}}{\partial \mathbf{u}}\left(\mathbf{u}_{i+\frac{1}{2}}^{n}\right) \right) \left(\frac{\partial {\mathbf{u}}}{\partial t}\right)_{i+\frac{1}{2}} = \widehat{\mathcal{A}}_{i+\frac{1}{2}} \left(\frac{\partial {\mathbf{u}}}{\partial t}\right)_{i+\frac{1}{2}}, 
\end{equation}
and the Generalized Riemann problem solver is used for the term $ \left(\frac{\partial {\mathbf{u}}}{\partial t}\right)_{i+\frac{1}{2}}$, given by 

\begin{equation}\label{eqn:dudt}
\left(\frac{\partial {\mathbf{u}}}{\partial t}\right)_{i+\frac{1}{2}} = - \widehat{\mathcal{A}}_{i+\frac{1}{2}}^{+} \Delta {\mathbf{u}}_{L, i+\frac{1}{2}} - \widehat{\mathcal{A}}_{i+\frac{1}{2}}^{-} \Delta {\mathbf{u}}_{R, i+\frac{1}{2}}, 
\end{equation}
where the operator $\Delta$ denotes the spacial derivative, such that $\Delta {\mathbf{u}}_{L/R, i+\frac{1}{2}}= \left(\frac{\partial {\mathbf{u}}}{\partial x}\right)_{L/R, i+\frac{1}{2}}$. The matrix $\widehat{\mathcal{A}}_{i+\frac{1}{2}}$ can be diagonized in the form of 
\begin{equation}\label{eqn:diago}
\widehat{\mathcal{A}}_{i+\frac{1}{2}} = \mathcal{R} \mathbf{\Lambda}\mathcal{L}, 
\end{equation}
where $\mathcal{R}$, and $\mathcal{L}$ denotes the right and left eigenmatrix, respectively, and $\mathbf{\Lambda}$ is the diagonal matrix consisting of eigenvalues $\lambda_k$. Then two unknown matrices in Eq.~\eqref{eqn:dudt} are defined as 
\begin{equation}\label{eqn:pm}
\widehat{\mathcal{A}}_{i+\frac{1}{2}}^{\pm} = \mathcal{R} \mathbf{\Lambda}^{\pm}\mathcal{L}, 
\end{equation}
where $\mathbf{\Lambda}^+$ consists of eigenvalues $\lambda^+_k = |\lambda_k| $, and  $\mathbf{\Lambda}^-$ of eigenvalues $\lambda^-_k = -|\lambda_k|$.

As aforementioned, we will focus on the high-order interpolation of flow variables $\mathbf{u}_{L/R,i+\frac{1}{2}}$ in section~\ref{sec:hermite}, and its first spacial derivative $\Delta {\mathbf{u}}_{L/R, i+\frac{1}{2}}$ in section~\ref{sec:hermite2}.

\subsection{Nonlinear interpolation of flow variables based on the Hermite polynomial}\label{sec:hermite}

In the proposed work we employ the Hermite polynomial instead of the Lagrangian approach for the higher-order approximation of the flow variables at the mid-point. The Hermite polynomial enjoys the advantage of its compact stencil, since it has been demonstrated that only half of the grid
points are required to derive the polynomial of the same degree when compared to the Lagrangian approach.

In the following work, we mainly consider the evaluation of variables on the left-hand side of $x_{i+\frac{1}{2}}$, i.e., ${\mathbf{u}}_{L,i+\frac{1}{2}}^n$ for the purpose of simplicity, and ${\mathbf{u}}_{R,i+\frac{1}{2}}^n$ is readily obtained from a symmetrical form of ${\mathbf{u}}_{L,i+\frac{1}{2}}$. We further drop the subscript $L$ and the supscript $n$ for the rest of the discussion, unless noted otherwise. 
%

As shown in Fig.~\ref{fig:beta}, given the node value ${\mathbf{u}}_{i}$, and its spatial derivative $\Delta {\mathbf{u}}_{i}$, only a three-point full stencil, which is denoted by $\text{S}_{i+\frac{1}{2}}=\{x_{i-1},x_{i},x_{i+1}\}$, is necessary for the construction of a fourth-degree Hermite polynomial $\mathbf{p}(x)$, such that $\mathbf{u}_{i+\frac{1}{2}} =\mathbf{p}\left(x_{i+\frac{1}{2}}\right)$. It can be determined by satisfying
\begin{equation}\label{eqn:full}
	\begin{array}{ll}
	\mathbf{p}\left(x_{k}\right)=\mathbf{u}_{k}, \; \; k=i-1, i, i+1 , \\
	\end{array}
\end{equation}

Then we can obtain a fifth-order spacial accurate interpolation for $\mathbf{u}_{i+\frac{1}{2}}$, which takes the form of 
\begin{equation}\label{eqn:fifth}
	\mathbf{u}_{i+\frac{1}{2}} = \mathbf{p}\left(x_{i+\frac{1}{2}}\right)=-\frac{1}{8} \mathbf{u}_{i-1}+\frac{9}{16} \mathbf{u}_{i}+\frac{9}{16} \mathbf{u}_{i+1}-\frac{3h}{64} \left( \Delta\mathbf{u}_{i-1}+3  \Delta\mathbf{u}_{i+1}\right).
\end{equation}

On the other hand, the full stencil can be split into three sub-stencils, namely, $\text{S}_{i+\frac{1}{2},1}=\left\{x_{i-1}, x_{i}\right\}, \text{S}_{i+\frac{1}{2},2}=\left\{x_{i}, x_{i+1}\right\},$ and $ \text{S}_{i+\frac{1}{2},3}=\left\{x_{i-1}, x_{i}, x_{i+1}\right\}$. Over each sub-stencil a specific Hermite polynomial can be constructed, which is subject to the following condition

\begin{equation}\label{eqn:interpo2}
	\begin{array}{ll}
		\mathbf{p}_{{1}}\left(x_{k}\right)=\mathbf{u}_{k}, & k=i-1, i, \quad \Delta\mathbf{p}_{{1}}\left(x_{i-1}\right)=\Delta \mathbf{u}_{i-1} \\
		\mathbf{p}_{{2}}\left(x_{k}\right)=\mathbf{u}_{k}, & k=i, i+1, \quad \Delta \mathbf{p}_{{2}}\left(x_{i+1}\right)=\Delta \mathbf{u}_{i+1} \\
		\mathbf{p}_{{3}}\left(x_{k}\right)=\mathbf{u}_{k}, & k=i-1, i, i+1
	\end{array}, 
\end{equation}
respectively. Once these three quadratic polynomials are determined, their corresponding interpolated mid-point flow variables are obtained as

\begin{equation}\label{eqn:poly}
	\begin{aligned}
		\mathbf{u}_{i+\frac{1}{2},{1}} &=\mathbf{p}_1\left(x_{i+\frac{1}{2}}\right)=-\frac{5}{4} \mathbf{u}_{i-1}+\frac{9}{4} \mathbf{u}_{i}-\frac{3h}{4} \Delta \mathbf{u}_{i-1} \\
		\mathbf{u}_{i+\frac{1}{2},{2}} &=\mathbf{p}_2\left(x_{i+\frac{1}{2}}\right)=\frac{1}{4} \mathbf{u}_{i}+\frac{3}{4} \mathbf{u}_{i+1}-\frac{h}{4} \Delta \mathbf{u}_{i+1} \\
		\mathbf{u}_{i+\frac{1}{2},{3}} &=\mathbf{p}_3\left(x_{i+\frac{1}{2}}\right)=-\frac{1}{8} \mathbf{u}_{i-1}+\frac{3}{4} \mathbf{u}_{i}+\frac{3}{8} \mathbf{u}_{i+1} 
	\end{aligned}. 
\end{equation}

Linear combination of $\mathbf{u}_{i+\frac{1}{2},{k}}$ yields 
\begin{equation}\label{eqn:linearcom}
	\mathbf{u}_{i+\frac{1}{2}} = \sum_{k=1}^3 d_k\mathbf{u}_{i+\frac{1}{2},{k}}, 
\end{equation}
where ${d}_k$s denote the linear weight, having values $d_{1}={1}/{16}, d_{2}={9}/{16}$, and $d_{3}={3}/{8} $, respectively. 

The interpolations in Eq.~\eqref{eqn:poly} can also be expressed in a generic form using (approximated) $n^{\,th}$ derivatives ($n=1,2$) in space, which are derived from the Taylor series, given by
\begin{equation}
	\mathbf{u}_{i+\frac{1}{2},k}=\mathbf{u}_i+ \mathbf{u}^{(1)}_{i,k} \frac{h}{2}+\mathbf{u}^{(2)}_{i,k}\frac{ h^2}{8}, \; \; k=1, 2, 3 ,
\end{equation}
where the first- and second-order derivatives are approximated by
\begin{equation} \label{eq:indicator}
	\begin{split}
		& \mathbf{u}^{(1)}_{i,0}=\frac{1}{h}\left(-2 \mathbf{u}_{i-1}+2 \mathbf{u}_{i}-h\Delta \mathbf{u}_{i-1}\right),\\
		& \mathbf{u}^{(1)}_{i,1}=\frac{1}{h}\left(-2 \mathbf{u}_{i}+2 \mathbf{u}_{i+1}-h\Delta \mathbf{u}_{i+1}\right),\\
		& \mathbf{u}^{(1)}_{i,2}=\frac{1}{2h}\left(-\mathbf{u}_{i-1}+\mathbf{u}_{i+1}\right),
	\end{split}
\end{equation}
\noindent and
\begin{equation} \label{eq:indicator2}
	\begin{split}
		& \mathbf{u}^{(2)}_{i,0}=\frac{1}{h^2}\left(-2\mathbf{u}_{i-1}+2\mathbf{u}_{i}-2h\Delta \mathbf{u}_{i-1}\right),\\
		& \mathbf{u}^{(2)}_{i,1}=\frac{1}{h^2}\left(2\mathbf{u}_{i}-2\mathbf{u}_{i+1}+2h\Delta \mathbf{u}_{i+1}\right),\\
		& \mathbf{u}^{(2)}_{i,2}=\frac{1}{h^2}\left(\mathbf{u}_{i-1}-2 \mathbf{u}_{i}+\mathbf{u}_{i+1}\right),
	\end{split}
\end{equation}
respectively. The smooth indicator, $\mathbf{\beta}_k$, is defined as \cite{Deng2000}
\begin{equation}
\beta_k=\left(h\mathbf{u}^{(1)}_{i,k}\right)^2+\left(h^2\mathbf{u}^{(2)}_{i,k}\right)^2, \!\!  \quad k= 1, 2, 3 \; .
\end{equation}

Nonlinear weight, $\omega_k$, is used to replace the linear weight, ${d}_k$, in Eq.~\eqref{eqn:linearcom} in order to alleviate non-physical oscillations when any sub-stencil is deemed crossed by a discontinuity. For instance, the nonlinear weight of Jiang and Shu~\cite{Jiang1996} {can be} used
\begin{equation}\label{eq:JS}
\omega_k=\frac{\alpha_k}{\sum_{k=0}^2\alpha_k}, \quad \alpha_k=\frac{d_k}{(\beta_k+\epsilon)^2},
\end{equation}
where the small parameter $\epsilon=10^{-6}$ is used to prevent division by zero.
It suggests that the corresponding JS weight can adaptively approach 0 for a substencil crossed by discontinuities, thus diminishing possible numerical oscillations, and continuously approximate the optimal linear weight in smooth regions, therefore achieving high-order accuracy.

\subsection{Nonlinear interpolation of the first spacial derivative of the flow variables}\label{sec:hermite2}
As shown in Eq.~\eqref{eqn:dudt}, the GRP solver requires the solution of the first spacial derivative in advance, which is the major focus of this section. Although a set of Hermite polynomials for the flow variables over the full stencil and the three sub-stencils have already been determined in the previous section, we do not recommend to simply differentiate them with respect to space and use the resulting solutions for the derivation of the first spacial derivatives involved; since it leads into one-order of accuracy loss in space, and even more concerning is that the convex combination strategy from the WENO concept cannot be used. A remedy is to construct a new set of Hermite polynomials for the computation of the first spacial derivative~\cite{Liu2015a,Liu2016}. First of all, a fifth-degree Hermite polynomial $\mathbf{p}(x)$ over the same three-point full stencil, $\text{S}_{i+\frac{1}{2}}=\{x_{i-1},x_{i},x_{i+1}\}$, is required to construct. It can be formulated by satisfying the particular conditions
\begin{equation}\label{eqn:full2}
\begin{array}{ll}
&\mathbf{p}\left(x_{k}\right)=\mathbf{u}_{k}, \; \; k=i-1, i, i+1 , \\
&\Delta \mathbf{p}\left(x_{l}\right)=\Delta \mathbf{u}_{l}, \; \; l=i-1, i, i+1. 
\end{array}
\end{equation}
Then a fifth-order spacial accurate interpolation for $\Delta \mathbf{u}_{i+\frac{1}{2}}$ can be obtained
\begin{equation}\label{eqn:fifth2}
\Delta \mathbf{u}_{i+\frac{1}{2}} =\Delta \mathbf{p}\left(x_{i+\frac{1}{2}}\right) =\frac{1}{h}\left( \frac{3}{64}\mathbf{u}_{i-1} -\frac{3}{2} \mathbf{u}_{i} + \frac{93}{64} \mathbf{u}_{i+1} \right) + \frac{1}{64}\left(\Delta\mathbf{u}_{i-1} -12\Delta\mathbf{u}_{i} - 15 \Delta\mathbf{u}_{i+1} \right).
\end{equation}

Over each of the three sub-stencils, namely, $\text{S}_{i+\frac{1}{2},1}=\left\{x_{i-1}, x_{i}\right\}, \text{S}_{i+\frac{1}{2},2}=\left\{x_{i}, x_{i+1}\right\},$ and $ \text{S}_{i+\frac{1}{2},3}=\left\{x_{i-1}, x_{i}, x_{i+1}\right\}$, a specific Hermite polynomial of third degree can be constructed, which is subject to the following condition
\begin{equation}\label{eqn:interpo22}
\begin{array}{ll}
\mathbf{p}_{{1}}\left(x_{k}\right)=\mathbf{u}_{k}, &  \Delta\mathbf{p}_{{1}}\left(x_{k}\right)=\Delta \mathbf{u}_{k}, \quad  k=i-1, i, \quad \\
\mathbf{p}_{{2}}\left(x_{k}\right)=\mathbf{u}_{k}, &  \Delta \mathbf{p}_{{2}}\left(x_{k}\right)=\Delta \mathbf{u}_{k}, \quad k=i, i+1, \quad \\
\mathbf{p}_{{3}}\left(x_{k}\right)=\mathbf{u}_{k}, & k=i-1, i, i+1, \quad \Delta \mathbf{p}_{{3}}\left(x_{i}\right)=\Delta \mathbf{u}_{i},
\end{array}
\end{equation}
respectively. The first spacial derivatives at the mid-point can be obtained thereafter

\begin{equation}\label{eqn:poly2}
\begin{aligned}
\Delta\mathbf{u}_{i+\frac{1}{2},{1}} &=\Delta\mathbf{p}_1\left(x_{i+\frac{1}{2}}\right)=\frac{9}{2h}\left( \mathbf{u}_{i-1} - \mathbf{u}_{i}\right) +\frac{7}{4} \Delta \mathbf{u}_{i-1}+\frac{15}{4} \Delta \mathbf{u}_{i}. \\
\Delta\mathbf{u}_{i+\frac{1}{2},{2}} &=\Delta\mathbf{p}_2\left(x_{i+\frac{1}{2}}\right)=\frac{3}{2h}\left( -\mathbf{u}_{i} + \mathbf{u}_{i+1}\right) -\frac{1}{4} \Delta \mathbf{u}_{i} -\frac{1}{4} \Delta \mathbf{u}_{i+1}.\\
\Delta\mathbf{u}_{i+\frac{1}{2},{3}} &=\Delta\mathbf{p}_3\left(x_{i+\frac{1}{2}}\right)=\frac{1}{8h}\left(  \mathbf{u}_{i-1}-8\mathbf{u}_{i} + 7\mathbf{u}_{i+1}\right) +\frac{1}{4} \Delta \mathbf{u}_{i}. 
\end{aligned}
\end{equation}

As is the case for $\mathbf{u}_{i+\frac{1}{2},{k}}$ in Eq.~\eqref{eqn:linearcom}, a convex combination can also be found for $\Delta\mathbf{u}_{i+\frac{1}{2},{k}}$, such that
\begin{equation}\label{eqn:linearcom2}
\Delta\mathbf{u}_{i+\frac{1}{2}} = \sum_{k=1}^3 \bar{d}_k \; \Delta\mathbf{u}_{i+\frac{1}{2},{k}}, 
\end{equation}
where the linear weight $\bar{d}_k$ is given by $\bar{d}_{1}={1}/{112}, \bar{d}_{2}={15}/{16}$, and $\bar{d}_{3}={3}/{56} $, respectively. 

The interpolations in Eq.~\eqref{eqn:poly} can also be expressed in a generic form using (approximated) $n^{\,th}$ derivatives ($n=2,3$),  given by
\begin{equation}
	\Delta \mathbf{u}_{i+\frac{1}{2},k}=\Delta \mathbf{u}_i+ \mathbf{u}^{(2)}_{i,k} \frac{h}{2}+\mathbf{u}^{(3)}_{i,k}\frac{ h^2}{8}, \; \; k=1, 2, 3 ,
\end{equation}
where the second- and third-order derivatives are approximated by
\begin{equation} \label{eq:indicator22}
\begin{split}
& \mathbf{u}^{(2)}_{i,0}=\frac{1}{h^2}\left(6\mathbf{u}_{i-1}-6\mathbf{u}_{i}+2h\Delta \mathbf{u}_{i-1} + 4h\Delta \mathbf{u}_{i}\right),\\
& \mathbf{u}^{(2)}_{i,1}=\frac{1}{h^2}\left(-6\mathbf{u}_{i}+6\mathbf{u}_{i+1}-4h\Delta \mathbf{u}_{i} - 2h\Delta \mathbf{u}_{i+1}\right),\\
& \mathbf{u}^{(2)}_{i,2}=\frac{1}{h^2}\left(\mathbf{u}_{i-1}-2 \mathbf{u}_{i}+\mathbf{u}_{i+1}\right),
\end{split}
\end{equation}
\noindent and
\begin{equation} \label{eq:indicator21}
\begin{split}
& \mathbf{u}^{(3)}_{i,0}=\frac{1}{h^3}\left(12 \mathbf{u}_{i-1}-12 \mathbf{u}_{i} + 6h\Delta \mathbf{u}_{i-1} + 6h\Delta \mathbf{u}_{i}\right),\\
& \mathbf{u}^{(3)}_{i,1}=\frac{1}{h^3}\left(12 \mathbf{u}_{i}-12 \mathbf{u}_{i+1} + 6h\Delta \mathbf{u}_{i} + 6h\Delta \mathbf{u}_{i+1}\right),\\
& \mathbf{u}^{(3)}_{i,2}=\frac{1}{h^3}\left(-3 \mathbf{u}_{i-1}+3 \mathbf{u}_{i+1} - 6h\Delta \mathbf{u}_{i}\right),
\end{split}
\end{equation}
respectively. The smooth indicator, $\mathbf{\beta}_k$, is defined as
\begin{equation}
\bar{\beta_k}=\left(h^2\mathbf{u}^{(2)}_{i,k}\right)^2+\left(h^3\mathbf{u}^{(3)}_{i,k}\right)^2, \!\!  \quad k= 1, 2, 3 \; .
\end{equation}

The nonlinear weight of Jiang and Shu~\cite{Jiang1996} {are also} used in this case, taking the form of 
\begin{equation}\label{eq:JS2}
\bar{\omega}_k=\frac{\alpha_k}{\sum_{k=0}^2\alpha_k}, \quad \alpha_k=\frac{d_k}{(\beta_k+\epsilon)^2},
\end{equation}
where the small parameter $\epsilon=10^{-6}$ is used to prevent division by zero.

\newpage
\begin{equation}\label{eqn:fifth-1st}
\Delta\mathbf{u}_{i+\frac{1}{2}} =\frac{1}{16h}\left( 3\mathbf{u}_{i-1} -24 \mathbf{u}_{i} + 21 \mathbf{u}_{i+1} + h\Delta\mathbf{u}_{i-1} -3h  \Delta\mathbf{u}_{i+1} \right), 
\end{equation}
for the the full stencil, and 
\begin{equation}\label{eqn:poly-1st}
\begin{aligned}
\Delta\mathbf{u}_{i+\frac{1}{2},{1}} &= \frac{1}{h}\left(-3\mathbf{u}_{i-1} + 3\mathbf{u}_{i}- 2h \Delta \mathbf{u}_{i-1}\right) \\
\Delta\mathbf{u}_{i+\frac{1}{2},{2}} &=\frac{1}{h}\left(-\mathbf{u}_{i}+  \mathbf{u}_{i+1}\right)\\
\Delta\mathbf{u}_{i+\frac{1}{2},{3}} &=\frac{1}{h}\left(-\mathbf{u}_{i}+  \mathbf{u}_{i+1}\right)
\end{aligned}. 
\end{equation}

%
\section{Numerical results} \label{sec:result}

\subsection{Sod and Lax shock tube problems}
Riemann initial-value problems of Sod~\cite{Sod1978} and Lax~\cite{Lax1954} are used to further evaluate shock-capturing capability of the proposed schemes employing discretization of the 1-D Euler equations
\begin{eqnarray*}
	\frac{\partial \rho}{\partial t} + \frac{\partial\big(\rho u \big)}{\partial x} \!\!\!\!&=&\!\!\! 0\, ,\\
	\frac{\partial\big(\rho u \big)}{\partial t} + \frac{\partial\big(\rho u^2\big)}{\partial x} \!\!\!\!&=&\!\!\! -\frac{\partial p}{\partial x}\, ,\\
	\frac{\partial E}{\partial t} + \frac{\partial\big(u E \big)}{\partial x} \!\!\!\!&=&\!\!\! -\frac{\partial\big(u\,p\big)}{\partial x}\, ,
\end{eqnarray*}
where $E = e + \frac{1}{2} u^2$ is the total energy per unit mass, and $e$ is internal energy.  The dependent variables are related through the perfect gas equation of state given by $p = (\gamma -1)\rho e$ with $\gamma = 1.4 $, thus closing the Euler equations system.

The Sod shock tube problem involves a right-moving shock of Mach number 1.7, while for the Lax shock tube problem, the right-moving shock has Mach number 2.0.  Initial conditions for the Sod problem are
\begin{equation}
	(\rho, u, p)=
	\begin{cases}
		\begin{matrix}
			(1,0,1)        &  \quad  x\in [0,0.5], \\
			(0.125,0,0.1)  & \quad  x\in  \left(0.5,1\right],
		\end{matrix}
	\end{cases}
\end{equation}
and the results at  $t=0.2$  are given by solving the problem on an evenly-distributed grid of $N=101$ points.

\begin{figure}[h!t]
	\begin{center}
		\subfigure[\label{fig:sod-density1}{}]{
			\resizebox*{8cm}{!}{\includegraphics{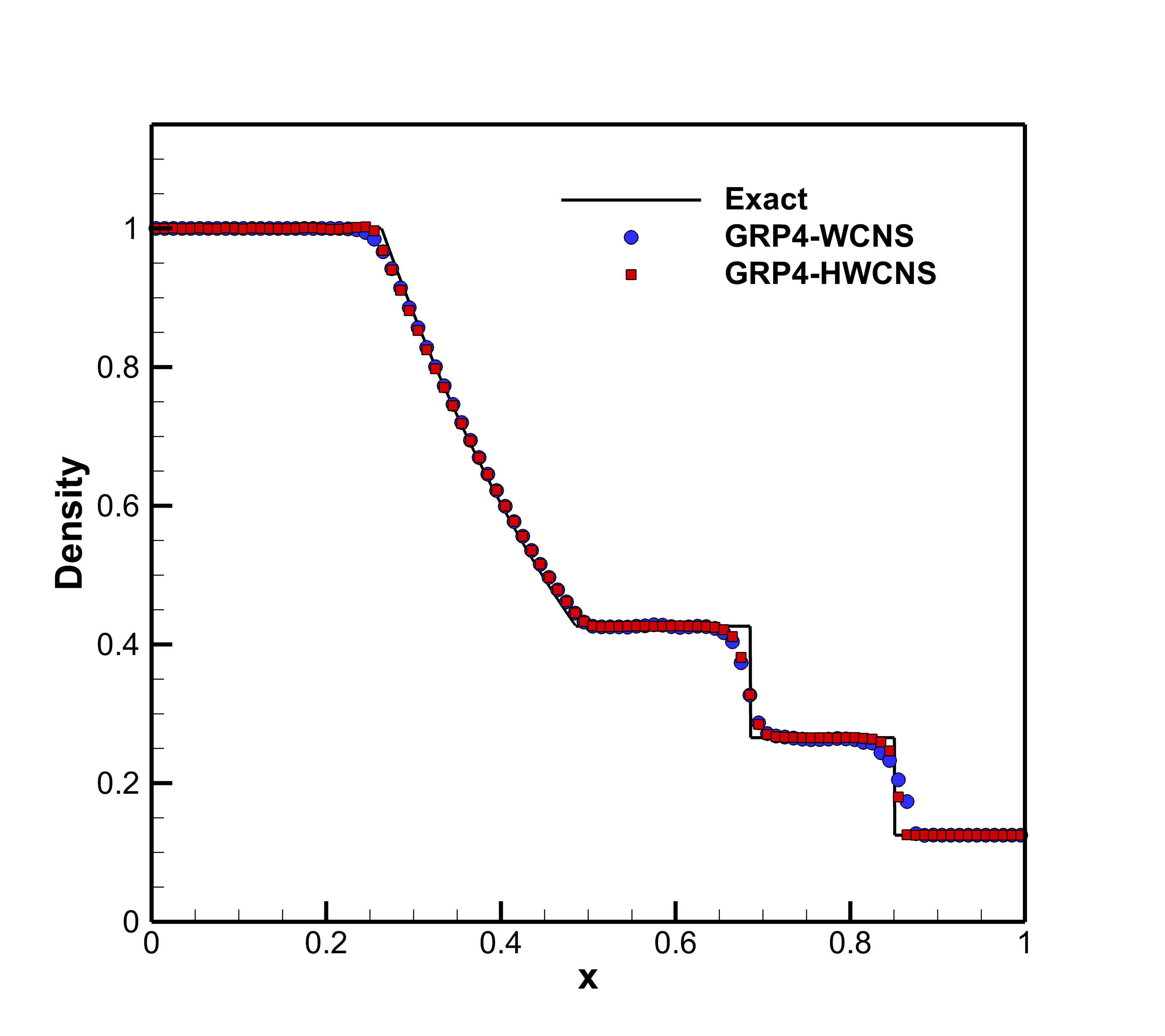}}}
		\subfigure[\label{fig:sod-density2}{}]{
			\resizebox*{8cm}{!}{\includegraphics{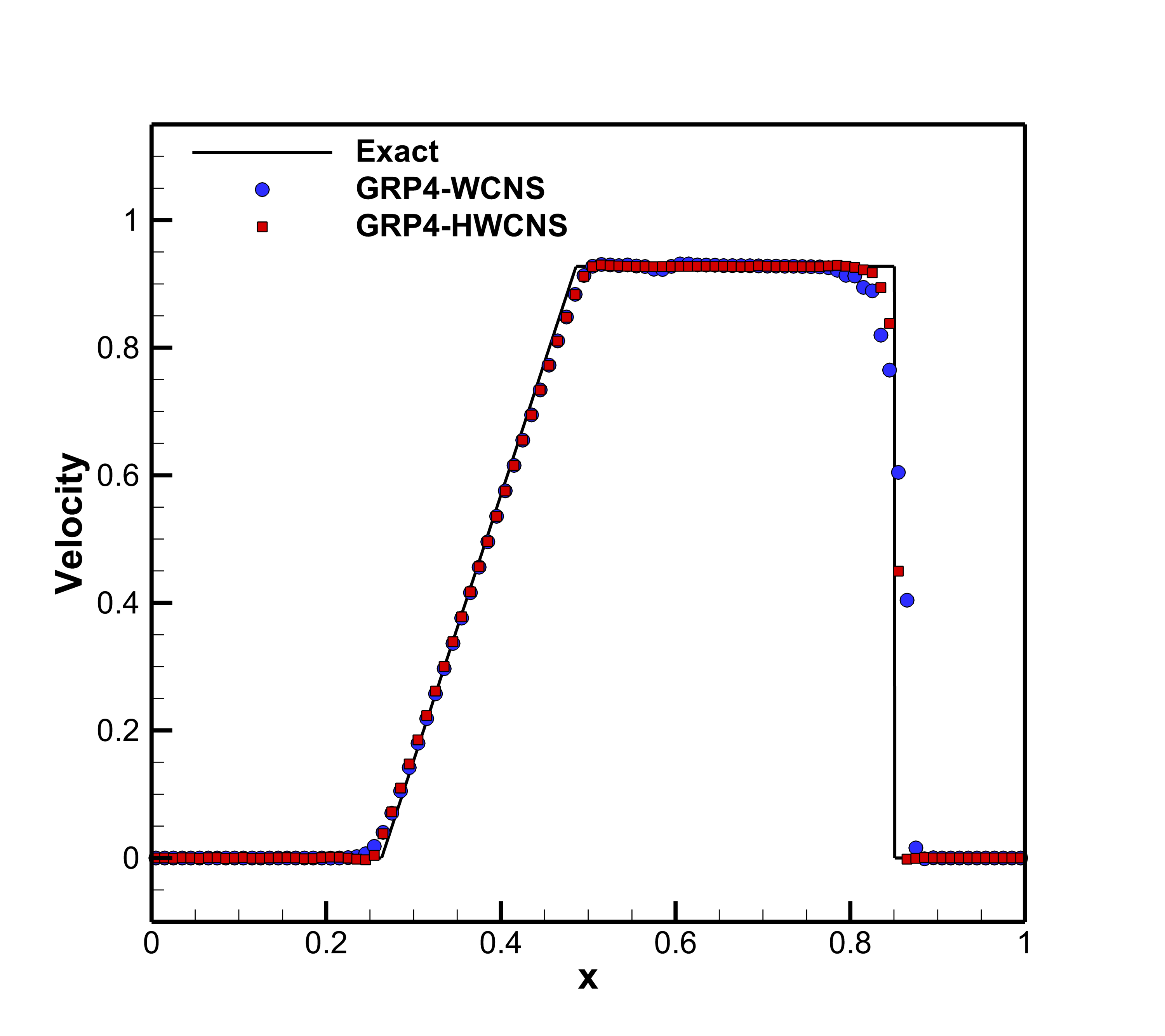}}}
		\caption{\label{fig:sod-density} Numerical and exact solutions for Sod problem at $t=0.2$; (a) density, (b) zoom in of density, (c) velocity.}
	\end{center}
\end{figure}

\begin{figure}[h!t]
	\begin{center}
		\subfigure[\label{fig:lax-density1}{\hspace{-8ex}}]{
			\resizebox*{8cm}{!}{\includegraphics{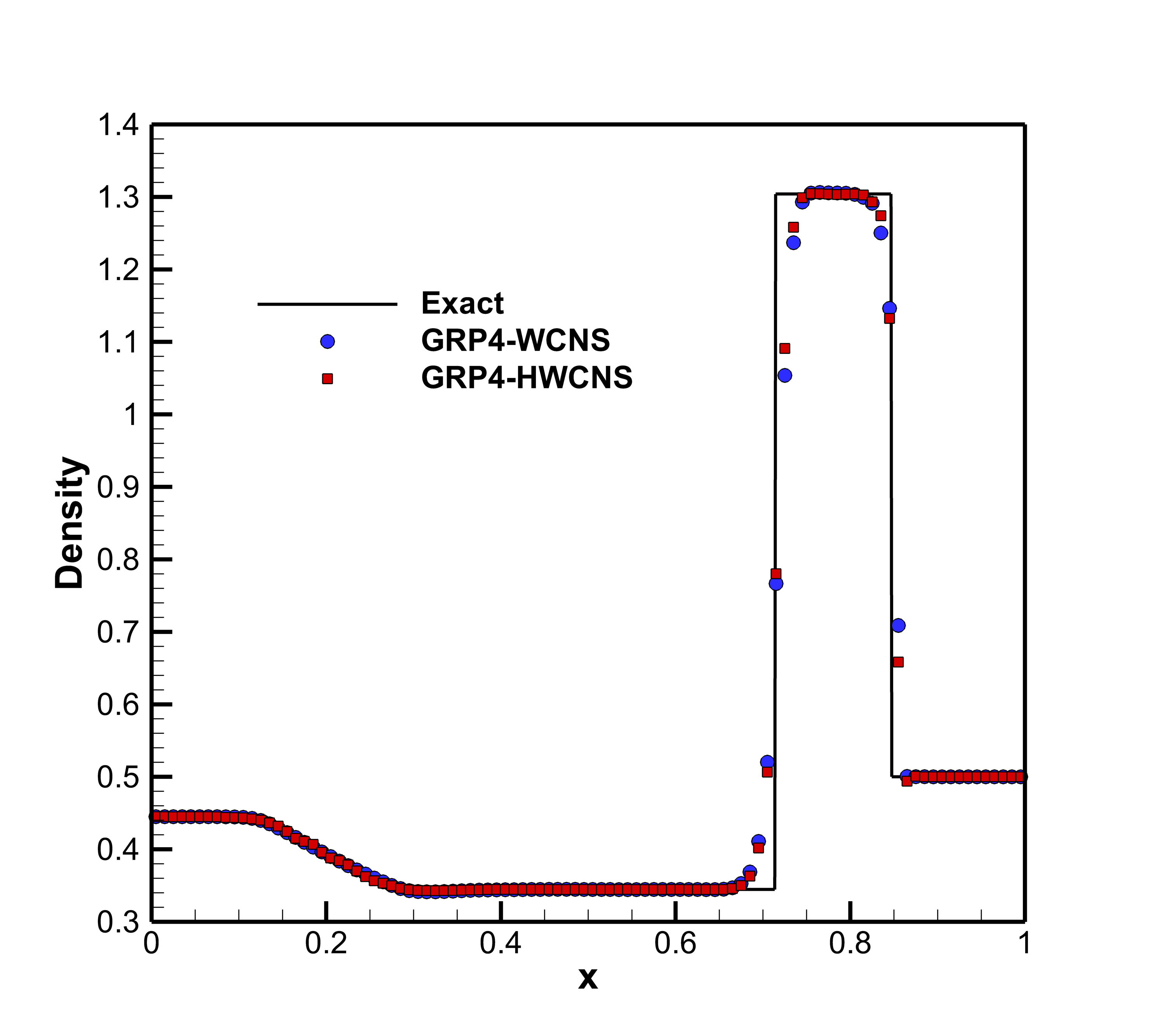}}}
		\subfigure[\label{fig:lax-density2}{\hspace{-8ex}}]{
			\resizebox*{8cm}{!}{\includegraphics{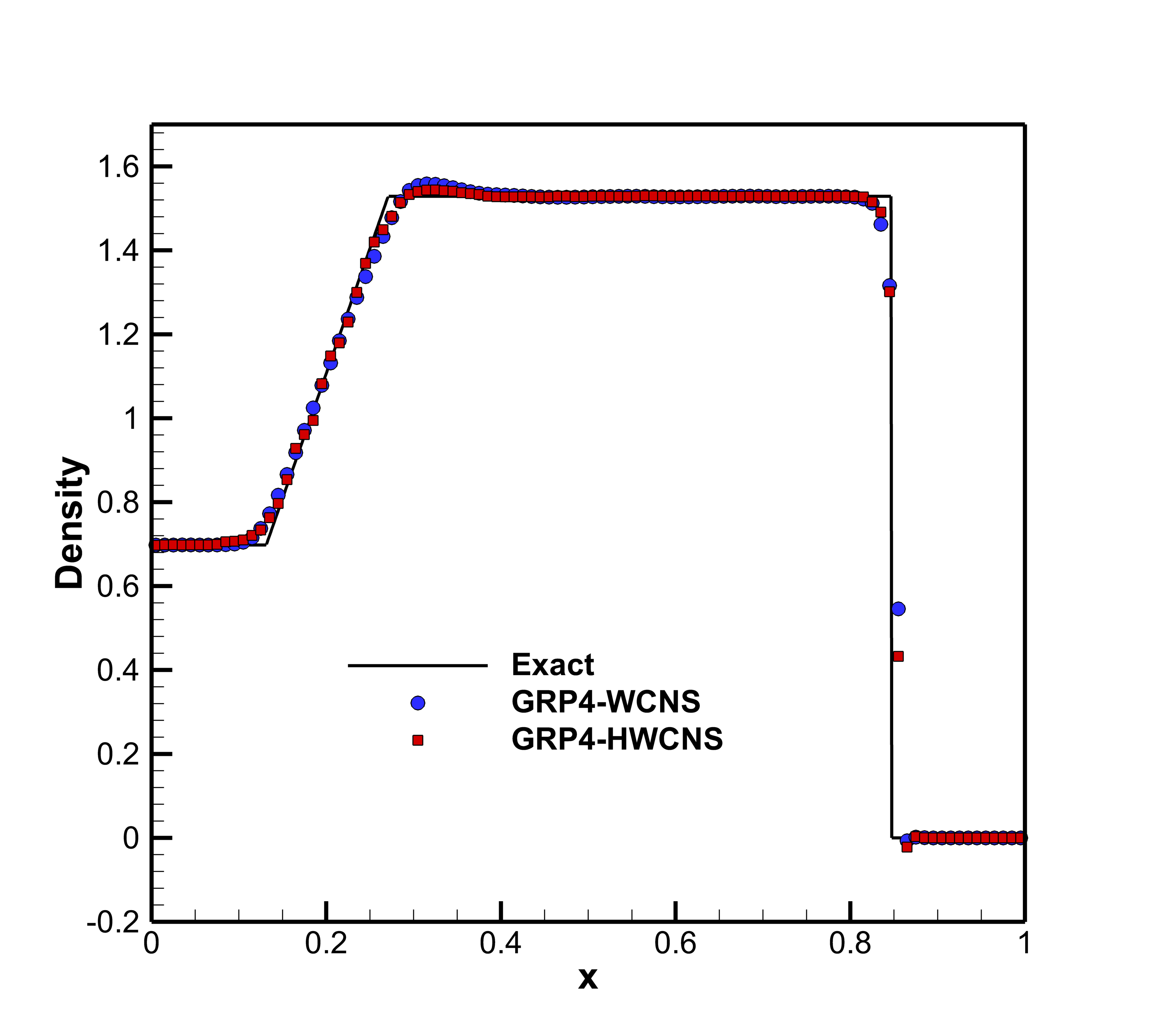}}}
		\caption{\label{fig:lax-density} Lax problem numerical and exact solutions at $t=0.14$; (a) density, and (b) velocity. }
	\end{center}
\end{figure}

Initial conditions for the Lax shock-tube problem are
\begin{equation}
	(\rho, u, p)=
	\begin{cases}
		\begin{matrix}
			(0.445, 0.698, 3.528)        &  \quad  x\in [0,0.5], \\
			(0.5,0,0.571)  & \quad  x\in  \left(0.5,1\right].
		\end{matrix}
	\end{cases}
\end{equation}
This case is also simulated on an evenly distributed grid of $N=101$ points, and the results at  $t=0.14$ are shown in Fig.\ \ref{fig:lax-density} for density and velocity distributions.

\subsection{Osher Shu}

\begin{figure}[H]
	\begin{center}
		\subfigure[\label{fig:shu-density1}{}]{
			\resizebox*{7.5cm}{!}{\includegraphics{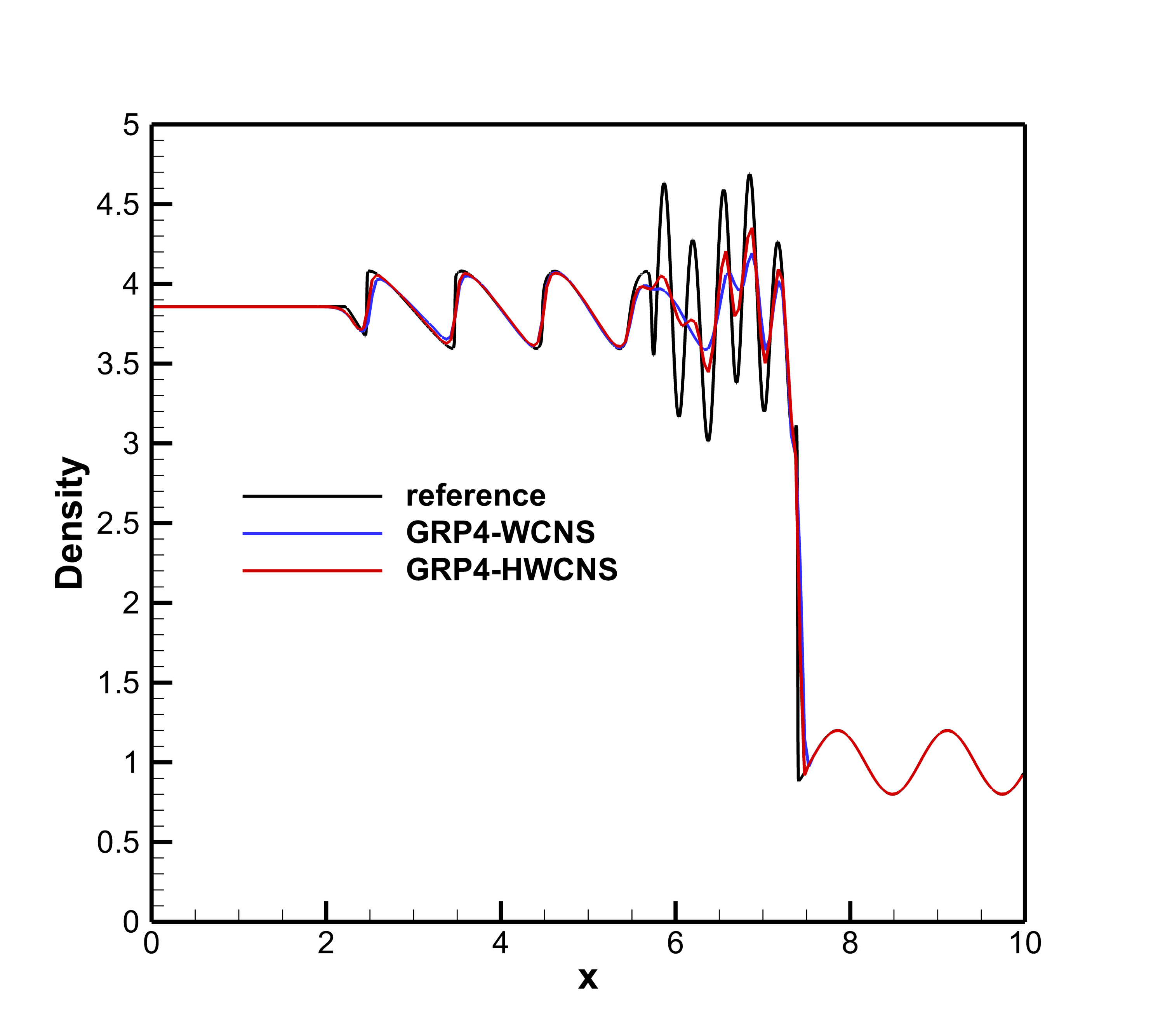}}}
		\subfigure[\label{fig:shu-density2}{}]{
			\resizebox*{7.5cm}{!}{\includegraphics{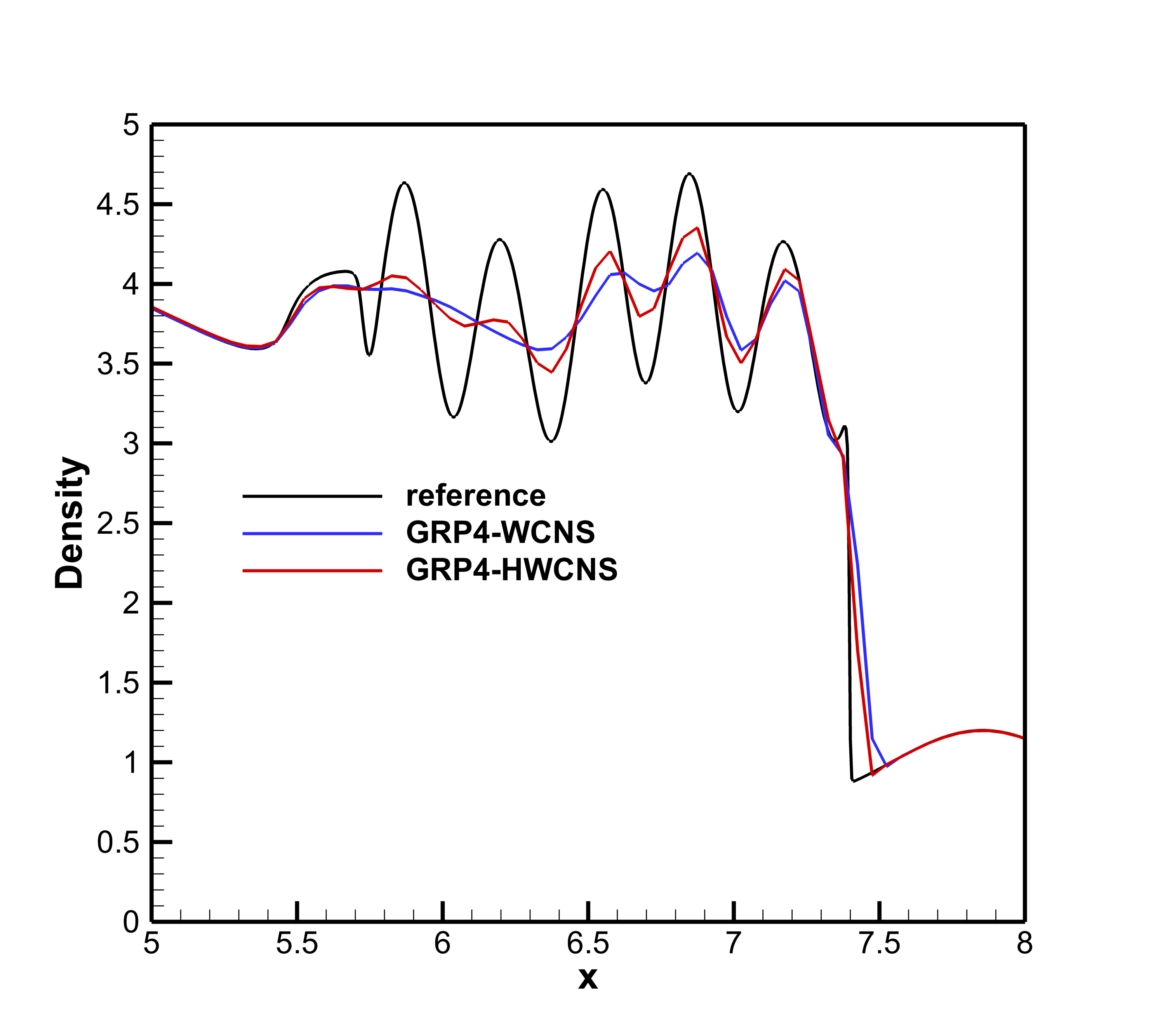}}}
		\caption{\label{fig:shu-density} Shock/density-wave interaction problem; numerical solutions and the exact solution at $t=5$; (a) full spatial domain, and (b) zoom in on high-amplitude region.}
	\end{center}
\end{figure}

\subsection{Two-dimensional Riemann problems}

\subsubsection{Configuration 6}

The equations of motion are the same as those used for configuration 3, the 2-D Euler equations, but now
with initial conditions given by
\begin{equation}
	(\rho, u, v, p)=
	\begin{cases}
		\begin{matrix}
			(1, 0.75, -0.5, 1)         &  \quad  (x,y)\in \left[\frac{1}{2},1\right]\!\times\!\left[\frac{1}{2},1\right], \\
			(2, 0.75, 0.5, 1)        & \quad   (x,y)\in \left[0 ,\frac{1}{2}\right)\!\times\!\left[\frac{1}{2}, 1\right], \\
			(1, -0.75, 0.5, 1)          &  \quad  (x,y)\in \left[0,\frac{1}{2}\right)\!\times\!\left[0,\frac{1}{2}\right), \\
			(3, -0.75, -0.5, 1)         & \quad   (x,y)\in \left[\frac{1}{2},1\right]\!\times\!\left[0,\frac{1}{2}\right). \\
		\end{matrix}
	\end{cases}
\end{equation}
Boundary conditions are the same as in the preceding test case.


\begin{figure}[H]
	\begin{center}
		\subfigure[\label{fig:rt-1}]{
			\resizebox*{8 cm}{!}{\includegraphics{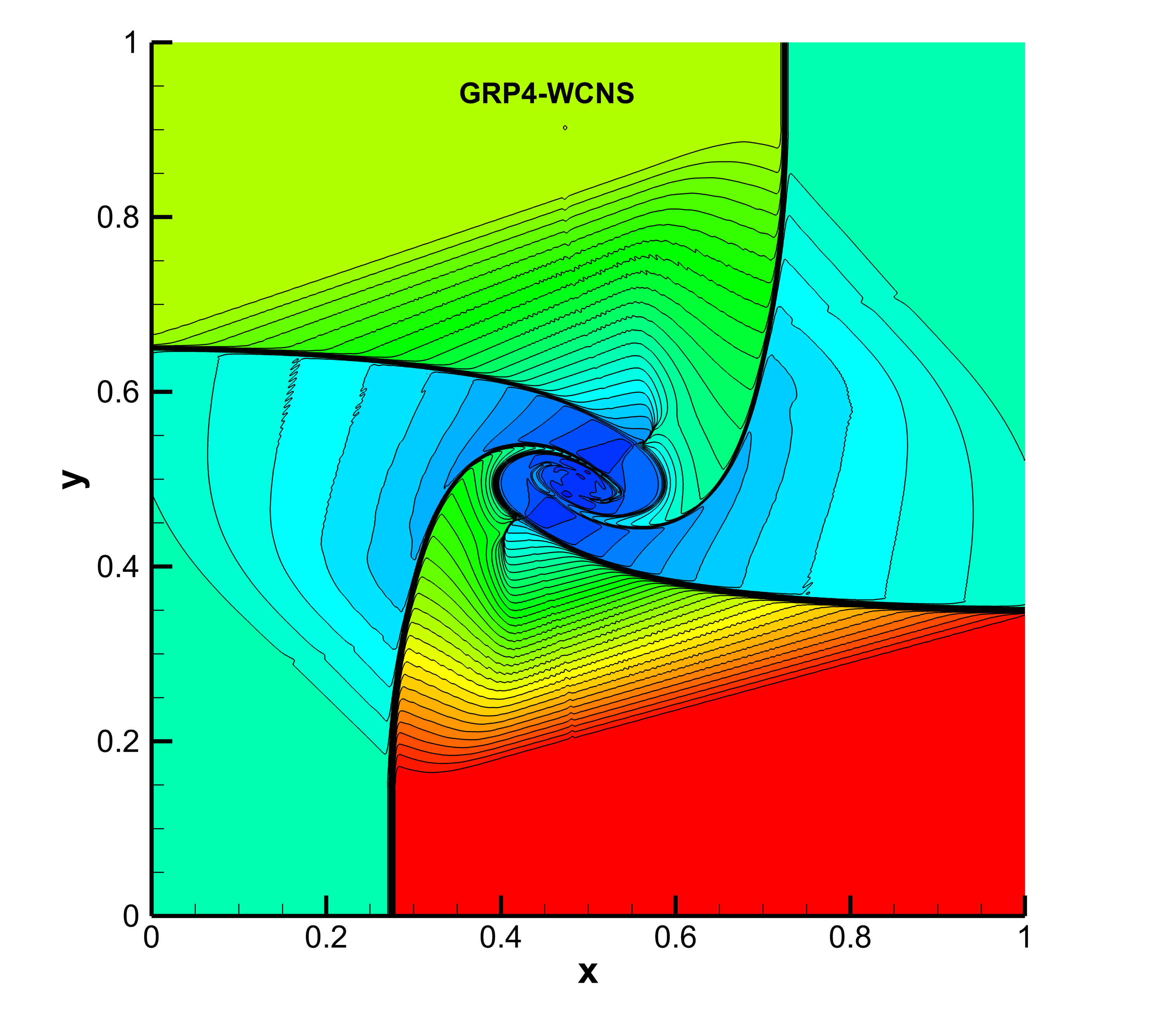}}}
		\subfigure[\label{fig:rt-2}]{
			\resizebox*{8 cm}{!}{\includegraphics{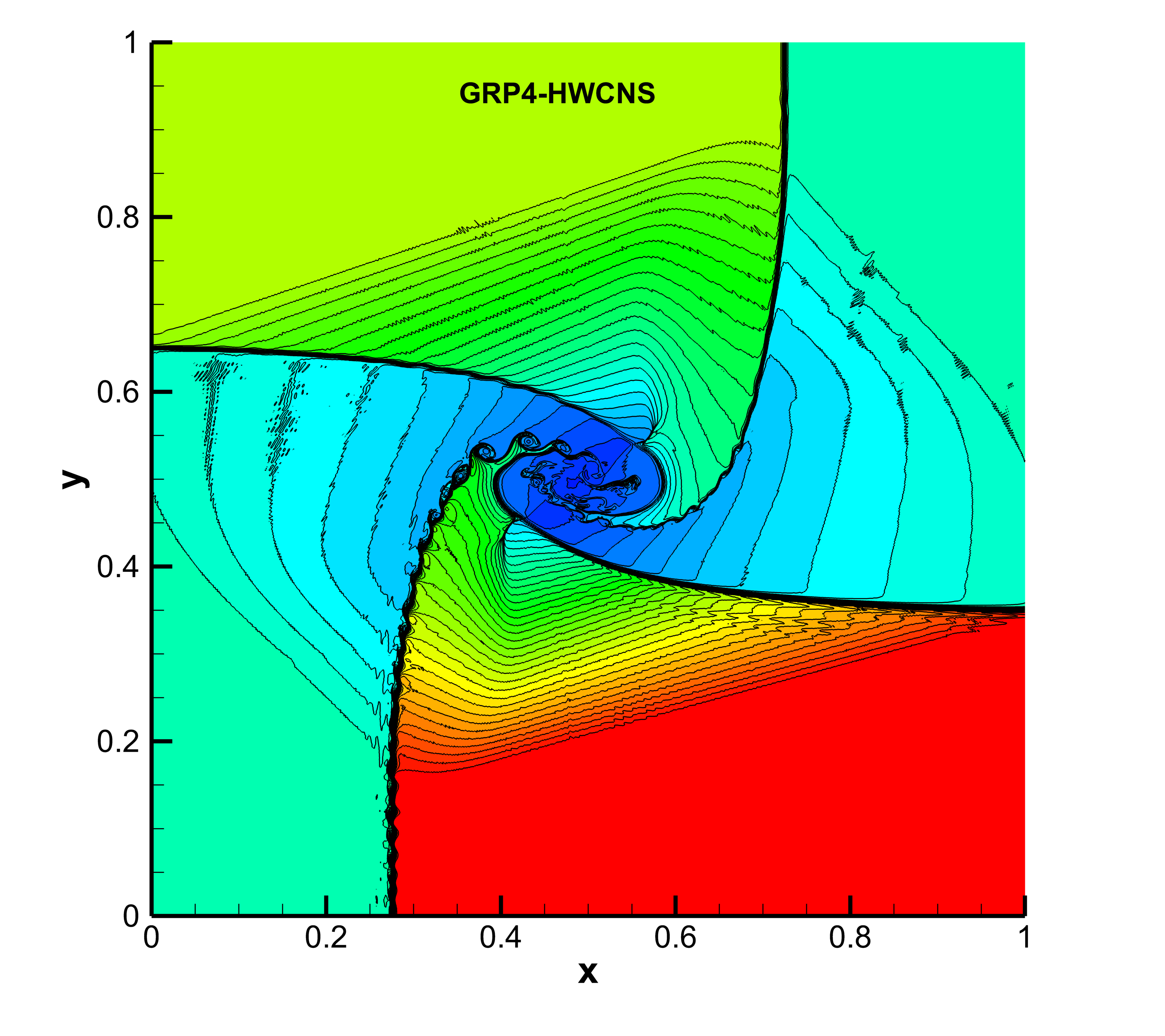}}}
	\caption{ \label{fig:2dR6} Configuration 6 of 2-D Riemann problems in \cite{Lax1998} computed on a grid of $1024\!\times\!1024$ points; 40 density contour lines ranging from 0.1 to 2.9 at $t=0.3$. }
	\end{center}
\end{figure}

\section*{Acknowledgments}
This work is supported by the open fund from the State Key Laboratory of Aerodynamics (Grant No.\ SKLA20180302), and the third author is partially supported by the National Natural Science Foundation of China (Grant No. 11872144).

\section{Conclusions} \label{sec:Conclusions}

\section*{Acknowledgments}
This work was financially supported by the National Numerical Wind Tunnel Project (No. NNW2019ZT5-B15) and the National Key Project (No. GJXM92579).

\bibliography{ref}

\end{document}